\documentclass[11pt,oneside,reqno]{amsart}

\usepackage{algorithm,algorithmicx,algpseudocode}
\usepackage{epsfig,psfrag,subfigure}
\usepackage{epstopdf}
\usepackage{amsmath,amssymb,amsfonts,latexsym, mathrsfs,euscript}
\usepackage{graphicx,graphics,cite}
\usepackage{wrapfig}
\usepackage[usenames]{color}
\usepackage{geometry}   
\usepackage{subfig}             
\usepackage{colortbl,booktabs,ctable,multirow}
\usepackage{epstopdf}
\usepackage{subfig}

\topmargin by -\headsep \textheight 8.7in \oddsidemargin 0.1in
\newcommand{\beq}{\begin{equation}}
\newcommand{\eeq}{\end{equation}}
\newcommand{\mc}{\multicolumn}
\newcommand{\cal}{\mathcal}

\newtheorem{theorem}{\bf Theorem}

\title[On Solving Ill-Conditioned Linear Systems]{On Solving Ill-Conditioned Linear Systems}

\author{Craig C. Douglas}
\address{School of Energy Resources, University of Wyoming, Laramie, WY 82071, USA}
  \email{craig.c.douglas@gmail.com}

\author{Long Lee}
\address{Department of Mathematics, University of Wyoming, Laramie, WY 82071-3036, USA}
\email{llee@uwyo.edu}

\author{Man-Chung Yeung}
\address{Department of Mathematics, University of Wyoming, Laramie, WY 82071-3036, USA}
\email{myeung@uwyo.edu}

\begin{document}

\maketitle


\begin{abstract}
This paper presents the first results to combine two theoretically sound
methods (spectral projection and multigrid methods) together to attack
ill-conditioned linear systems.
Our preliminary results show that the proposed algorithm applied to a Krylov
subspace method takes much fewer iterations for solving an ill-conditioned
problem downloaded from a popular online sparse matrix collection.
\end{abstract}


%
%


\section{Introduction}
\label{sect:introduction}

It is well-known that robustness and efficiency of iterative methods are
affected by the condition number of a linear system.
When a linear system has a large condition number, usually due to eigenvalues
that are close to the origin of the spectrum domain, iterative methods tend to
take many iterations before a convergence criterion is satisfied.
Sometimes, iterative methods will fail to converge within a reasonable
computer elapsed time, or even do not converge at all, if the condition number
is too large.
Unstable linear systems, or systems with large condition numbers, are called
ill-conditioned.
For an ill-conditioned linear system, slight changes in the coefficient matrix
or the right-hand-side cause large changes in the solution.
Typically, roundoff error in the computer arithmetics can cause instability
when attempts are made to solve an ill-conditioned system either directly or
iteratively on a computer.

It is widely recognized that linear systems resulting from discretizing
ill-posed integral equations of the first kind are highly ill-conditioned.
This is because the eigenvalues for the first kind integral equations with
continuous or weakly singular kernels have an accumulation point at zero.
Integral equations of the first kind are frequently seen in statistics, such
as unbiased estimation, estimating a prior distribution on a parameter given
the marginal distribution of the data and the likelihood, and similar tests
for normal theory problems.
They also arise from indirect measurements and nondestructive testing in
inverse problems. Other  ill-conditioned linear systems can be seen in
training of neural networks, seismic analysis, Cauchy problem for
parabolic equations, and multiphase flow of chemicals.
For pertinent references of ill-conditioned linear systems,
one can see, for example, Engl \cite{engl93} and Groetsch \cite{groetsch97}.


Solving these ill-conditioned linear algebra problems becomes a long-standing
bottleneck for advancing the use of iterative methods.
The convergence of iterative methods for ill-conditioned problems, however,
can be improved by using preconditioning.
Development of preconditioning techniques is therefore a very active research
area.
A preconditioning strategy that {\em deflates} few isolated external
eigenvalues was first introduced by Nicolaides \cite{nicolaides87}, and
investigated by several others \cite{mansfield91,tv06, vsm99,fv01}.
The deflation strategy is an action that removes the influence of a subspace
of the eigenspace on the iterative process.
A common way to deflate an eigenspace is to construct a proper projector $P$
as a preconditioner and solve
\begin{equation}\label{equ:1-28-1}
PAx= Pb,\quad P, A\in\mathbb{C}^{N\times N}.
\end{equation}
The deflation projector $P$, orthogonal
the matrix $A$ and the vector $b$
against some subspace, is defined by
\begin{equation}\label{eq:projector}
P=I-AZ(Z^{H}AZ)^{-1}Z^{H}, \quad Z\in\mathbb{C}^{N\times m},
\end{equation}
where $Z$ is a matrix of deflation subspace, i.e., the space to be projected out of
the residual, and $I$ is the identity matrix of appropriate size \cite{saad96,
fv01}.
We assume that (1) $m \ll N$ and (2) $Z$ has rank $m$.
A deflated $N{\times}N$ system (\ref{equ:1-28-1}) has an eigensystem different from that
of $Ax =b$.
Suppose that $A$ is diagonalizable, and set $Z= [v_1,\cdots,v_m]$, whose
columns are eigenvectors of $A$ associated with eigenvalues
$\lambda_1,\cdots,\lambda_m$.
Then the spectrum $\sigma(PA)$ would contain the same eigenvalues of $A$,
except $\lambda_1,\cdots,\lambda_m$.
Usually, eigenvectors are not easily available.
This motivates us to develop an efficient and robust algorithm for finding an
approximate deflation subspace, other than using the exact eigenvectors to
construct the deflation projector $P$.

Suppose that we want to deflate a set of eigenvalues of $A$ enclosed in a
circle $\Gamma$ that is centered at the origin with the radius $r$.
Without loss of generality, let this set of eigenvalues be
$\{\lambda_1,\cdots,\lambda_k\}$.
Let the subspace spanned by the
corresponding eigenvectors of $\{\lambda_1,\cdots,\lambda_k\}$ be
$\mathcal{Z}_k =\text{Span}\{v_1,\cdots,v_k\}$.
Then the deflation subspace matrix $Z$ in (\ref{eq:projector}) obtained
by randomly
selecting $m$
vectors from
$\mathcal{Z}_k$ can be written as a contour integral\cite{saad11}
\beq\label{eq:Z_oint}
Z=\frac{1}{2\pi \sqrt{-1}}\oint_{\Gamma}(zI-A)^{-1}Y\,dz,
\eeq
where $Y$ is a random matrix of size $N{\times}m$.
If the above contour integral is approximated by a Gaussian quadrature, we
have
\beq\label{eq:gaussian}
Z=\sum_{i=1}^{q}\omega_i(z_iI-A)^{-1}Y,
\eeq
where $\omega_i$ are the weights, $z_i$ are the Gaussian points, and $q$ is
the number of Gaussian points on $\Gamma$ for the quadrature.
It is worth noting that (\ref{eq:gaussian}) is required to solve $q$
shifted linear systems $(z_i I-A)X= Y$, $i=1,\cdots,q$.
Using (\ref{eq:gaussian}) for the deflation projector $P$ in
(\ref{eq:projector}), the preconditioned linear system (\ref{equ:1-28-1}) is no longer
severely ill-conditioned.

We remark that the construction of a deflation subspace matrix $Z$ through
(\ref{eq:gaussian}) is motivated by the works in \cite{ss03, st07, polizzi,
TP13}.


\section{Methodology}
\label{sec:problem}

We consider the solution of the linear system
\begin{equation}\label{equ:11-19-1}
A x = b
\end{equation}
by a Krylov subspace method, where we assume that $A \in {\mathbb C}^{N \times N}$ is
nonsingular and $b \in {\mathbb C}^N$.
Let an initial guess $x_0 \in {\mathbb C}^N$ be given along with its residual $r_0 = b - A x_0$. A Krylov subspace method
recursively constructs an approximate solution, $x_j$, such that
$$
x_j \in x_0 + {\cal K}_j(A, r_0) \equiv x_0 + \mbox{span} \{r_0, A r_0, \ldots, A^{j-1} r_0\},
$$
and its residual $r_j = b - Ax_j$ satisfies some desired conditions. It is well-known that the convergence rate of a Krylov subspace method depends on the eigenvalue distribution of the coefficient matrix $A$. A variety of error bounds on $r_j$ exist in the literature. Let us take GMRES\cite{Saa86S} as an example.


\subsection{GMRES}
\label{Subsec:Problem:GMRES}

In GMRES, the residual $r_j$ is required to satisfy the condition
$$
\|r_j \|_2 = \min_{\xi \in x_0 + {\cal K}_j(A, r_0)}\|b - A \xi\|_2,
$$
namely, the approximate solution $x_j$ obtained at iteration $j$ of GMRES is optimal in terms of residual norm.
In the case where $A$ is
diagonalizable, an upper bound on $\|
r_j\|_2$ is provided by the following result.

\begin{figure}[tbh]
\begin{center}
\includegraphics[width=4.8in]{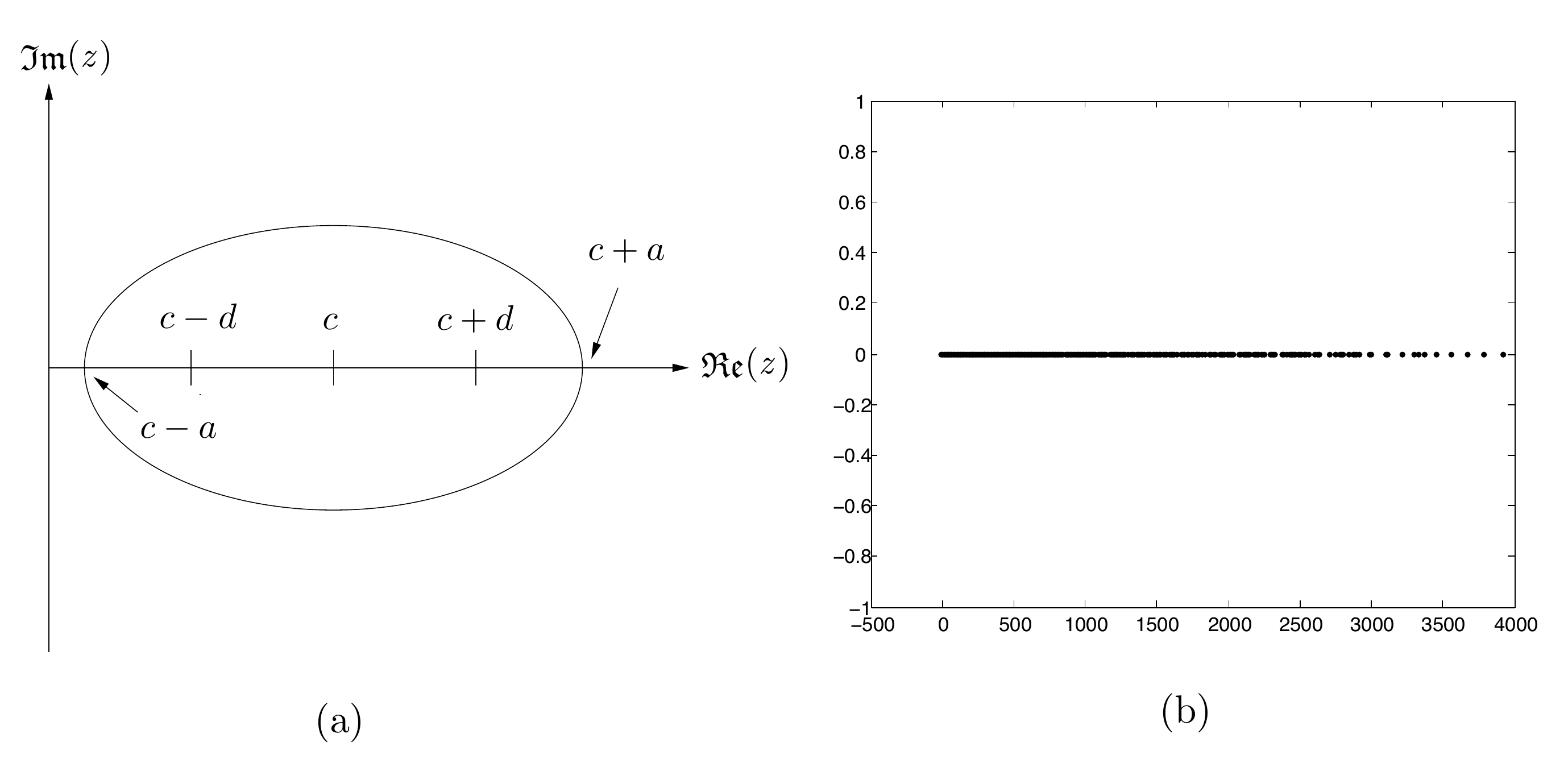}
\end{center}


\caption{(a) A schematic ellipse in the complex plane with center $c$,  focal distance $d$, and semi-major axis $a$. (b) Eigenvalue distribution of the test matrix {\it bcsstm27}.}
\label{fig:ellipse}
\end{figure}

\begin{theorem}\label{cor:3-16-1} (\cite[Corollary 6.33]{saad96})
Suppose that $A$ can be decomposed as
\begin{equation}\label{equ:11-19-2}
A = V \Lambda V^{-1}
\end{equation}
with $\Lambda$ being the diagonal matrix of eigenvalues.
Let
$E(c, d, a)$ denote the ellipse in the complex plane with center $c$,
focal distance $d$,
and semi-major axis $a$ (see Fig. \ref{fig:ellipse}(a)).
If all the eigenvalues of $A$ are located in $E(c, d, a)$ that excludes the origin of the complex plane, then
\begin{equation}\label{equ:3-16-31}
\|r_j\|_2 \leq \kappa_2(V) \frac{C_j(\frac{a}{d})}{|C_j(\frac{c}{d})|} \|r_0\|_2
\end{equation}
where $\kappa_2(V) = \|V\|_2 \|V^{-1}\|_2$ and $C_j$ is the Chebyshev polynomial of degree $j$.
\end{theorem}


An explicit expression of $C_j(\frac{a}{d})/C_j(\frac{c}{d})$ can be found on p.207 of \cite{saad96}, and under some additional assumptions on $E(c, d, a)$
(say, the ellipse in Fig. \ref{fig:ellipse}(a)
)
\begin{equation}\label{equ:3-18-153}
\frac{C_j(\frac{a}{d})}{C_j(\frac{c}{d})} \approx \left(\frac{a + \sqrt{a^2 - d^2}}{c + \sqrt{c^2 - d^2}}\right)^j \equiv \delta^j.
\end{equation}

The upper bound in (\ref{equ:3-16-31}) therefore contains two factors: the
condition number $\kappa_2(V)$ of the eigenvector matrix $V$
 and the scalar $\delta$ determined by the
distribution of the eigenvalues of $A$. If $A$ is nearly normal and has a
spectrum $\sigma(A)$ which is clustered around $1$, we would have
$\kappa_2(V) \approx 1$ and $\delta < 1$. In this case,
$\|r_j\|_2$ decays exponentially in a rate of power $\delta^j$,
resulting in a fast convergence of GMRES. The error bound (\ref{equ:3-16-31}) does not tell the whole story, however, because the convergence rate can also be enhanced if the eigenvalues of $A$ are clustered\cite{vdS86v}.


Since the ellipse $E(c, d, a)$ in Theorem~\ref{cor:3-16-1} is
required to include all eigenvalues of $A$, the outlying
eigenvalues may keep the ellipse large, implying a large $\delta$.
To reduce $\delta$,
we therefore wish to remove these outlying
eigenvalues from  $\sigma(A)$. Any procedure of
doing so is known as {\it deflation}. GMRES in combination with deflation is called {\it Deflated GMRES}. 


\subsection{Deflated GMRES}
\label{sec:conv}

Suppose $x^*$ is the exact
solution of (\ref{equ:11-19-1}). Let a so-called deflation-subspace matrix $Z = [z_1, \ldots, z_m] \in {\mathbb
C}^{N \times m}$ be given, whose columns are linearly independent.
Define
the two projectors\cite{vsm99, fv01}
\begin{equation} \label{equ:8-21-1}
\begin{array}{rcl} P \equiv I - A Z (Z^H A Z)^{-1} Z^H &
\mbox{and} & \widetilde{P} \equiv I - Z (Z^H A Z)^{-1} Z^H A,
\end{array}
\end{equation}
where $Z^H A Z$ is assumed to be invertible. It is straightforward to
verify that $P^2 = P, \; \widetilde{P}^2 = \widetilde{P}$ and $P A = A
\widetilde{P}$.

Using $\widetilde{P}$, we split $x^*$ into two parts:
$$
x^* = (I - \widetilde{P}) x^* + \widetilde{P} x^* \equiv x^*_1 + x^*_2.
$$
For $x^*_1$, we have
$$
x^*_1 = (I - \widetilde{P}) x^* = Z (Z^H A Z)^{-1} Z^H A x^* = Z (Z^H A
Z)^{-1} Z^H b.
$$
For
$x^*_2$, we obtain
$$
x^*_2 = A^{-1} P b,
$$
since
$
A x^*_2 = A \widetilde{P} x^* = P A x^* = P b.
$
Now, if $x^\#$ is a solution of the singular system
\begin{equation}\label{equ:7-31}
P A x = P b,
\end{equation}
then
$$\begin{array}{rcl}
A \widetilde{P} x^\# = P b & \Leftrightarrow & \widetilde{P} x^\# = A^{-1} P b
 = x^*_2.
\end{array}
$$
Based on the above observation, a Deflated GMRES
algorithm is given in Algorithm~\ref{alg:DGMRES}.

\begin{algorithm}
\begin{algorithmic}
\State Choose $Z$; \\
\State Compute $x^*_1 = Z(Z^H A Z)^{-1} Z^H b$; \\
\State Solve $P A x = P b$ by GMRES to obtain a solution $x^\#$; \\
\State Compute $x_2^* = \widetilde{P} x^\#$; \\
\State Determine $x^* = x^*_1 + x^*_2$.
\end{algorithmic}
\caption{Deflated GMRES}
\label{alg:DGMRES}
\end{algorithm}

Assume that the nonsingular
 $A \in {\mathbb C}^{N \times N}$ has a
decomposition~(\ref{equ:11-19-2}) with $V = [v_1, \ldots, v_N]$ and $\Lambda = \mbox{diag}\{\lambda_1, \ldots, \lambda_N\}$.
If we set $Z = [v_1, \ldots, v_m]$ in (\ref{equ:8-21-1}),
then the spectrum $\sigma(PA)$ contains all the eigenvalues of $A$ except $\lambda_1, \ldots, \lambda_m$, namely, $\sigma(PA) = \{0, \cdots, 0, \lambda_{m+1}, \cdots, \lambda_N\}$.

Perform a $QR$ factorization on $V$ as follows:
\begin{equation}\label{equ:3-17-101}
V = Q R \equiv [Q_1, Q_2] \left[
\begin{array}{cc}
R_{11} & R_{12}\\
0 & R_{22}
\end{array} \right],
\end{equation}
where $Q_1 \in {\mathbb C}^{N \times m}$ and $R_{11} \in {\mathbb C}^{m \times m}$.
If we set $Z = [v_1, \ldots, v_m]$ and apply GMRES to solve (\ref{equ:7-31}), an upper bound on $\|r_j\|$ is given by the following theorem\cite{YTV}.

\begin{theorem}\label{cor:4-4-50} Suppose that $A$ has a decomposition (\ref{equ:11-19-2}), and
suppose GMRES is used to solve (\ref{equ:7-31}) with $Z = [v_1, \ldots, v_m]$. 
If all the eigenvalues $\lambda_{m+1}, \ldots, \lambda_N$ of $A$ are located in an ellipse $E(c, d, a)$ which excludes the origin of the complex plane, then
\begin{equation}\label{equ:3-18-32}
\|r_j\|_2 \leq \kappa_2(R_{22}) \frac{C_j(\frac{a}{d})}{|C_j(\frac{c}{d})|} \|r_0\|_2. 
\end{equation}
\end{theorem}

With (\ref{equ:3-18-153}),
the upper bound (\ref{equ:3-18-32}) of the residual norm $\|r_j\|_2$ of Deflated GMRES is determined by the condition number of $R_{22}$ (rather than $V$), and the scalar $\delta$ which is determined by the distribution of the undeflated eigenvalues $\lambda_{m+1}, \ldots, \lambda_N$ of $A$.


\subsection{Spectral Projector and Construction of $Z$}
\label{Subsec:Problem:SpectralProj}

Spectral projector is described in detail in \S3.1.3-\S3.1.4 of \cite{saad11}. Other reference includes \cite{chatelin, dunford, kato}.
Let $A = VJV^{-1}$
be the Jordan canonical decomposition of $A$ where
$$\begin{array}{ccc}
V = [v_1, v_2, \ldots, v_N] & \mbox{and} & J = diag\{J_{N_1}(\lambda_1), J_{N_2}(\lambda_2), \ldots, J_{N_d}(\lambda_d)\}.
\end{array}
$$
The eigenvalues $\lambda_i$ in $J$ are not necessarily distinct and can be repeated according to their multiplicities, and the diagonal block $J_{N_i}(\lambda_i)$ in $J$ is an $N_i \times N_i$ Jordan block associated with the eigenvalue $\lambda_i$.

Let $\Gamma$ be a given positively oriented simple closed curve in the complex plane. Without loss of generality, let the set of eigenvalues of $A$ enclosed by $\Gamma$ be $\{\lambda_1, \lambda_2, \ldots, \lambda_k\}$. In other words, the eigenvalues $\lambda_{k+1}, \ldots, \lambda_d$ lie outside the region enclosed by $\Gamma$. Set $s \equiv N_1 + N_2 + \ldots + N_k$, the number of eigenvalues inside $\Gamma$ with multiplicity taken into account. Then the residue
$$
P_\Gamma = \frac{1}{2 \pi \sqrt{-1}} \oint_\Gamma (zI - A)^{-1} dz
$$
defines a projection operator onto the space $\sum_{i = 1}^k \mbox{Null}(A - \lambda_i I)^{l_i}$ where $l_i$ is the index of $\lambda_i$, namely,
$$
\mbox{Range}(P_\Gamma) = \mbox{span}\{v_1, v_2, \ldots, v_s\}.
$$
In particular, if $A$ has a diagonal decomposition (\ref{equ:11-19-2}), $P_\Gamma$ is a projector onto the sum $\sum_{i = 1}^k {\mathbb E}_{\lambda_i}$ of the $\lambda_i$-eigenspace ${\mathbb E}_{\lambda_i}$ of $A$.

Pick a random matrix $Y \in {\mathbb C}^{N \times s}$ and set
\begin{equation}\label{equ:11-24-1}
Z = P_\Gamma Y = \frac{1}{2 \pi \sqrt{-1}} \oint_\Gamma (zI - A)^{-1} Y dz
\end{equation}
in (\ref{equ:8-21-1}). Then we almost surely have $\sigma(PA) = \{0, \cdots, 0, \lambda_{k+1}, \cdots, \lambda_d\}$. Therefore all the eigenvalues of $A$ inside $\Gamma$ are removed from the spectrum of $PA$.


\subsection{Numerical Examples}
\label{Subsec:Problem:NumerExI}

\begin{figure}[tbh]
\begin{center}
\includegraphics[width=4.85in]{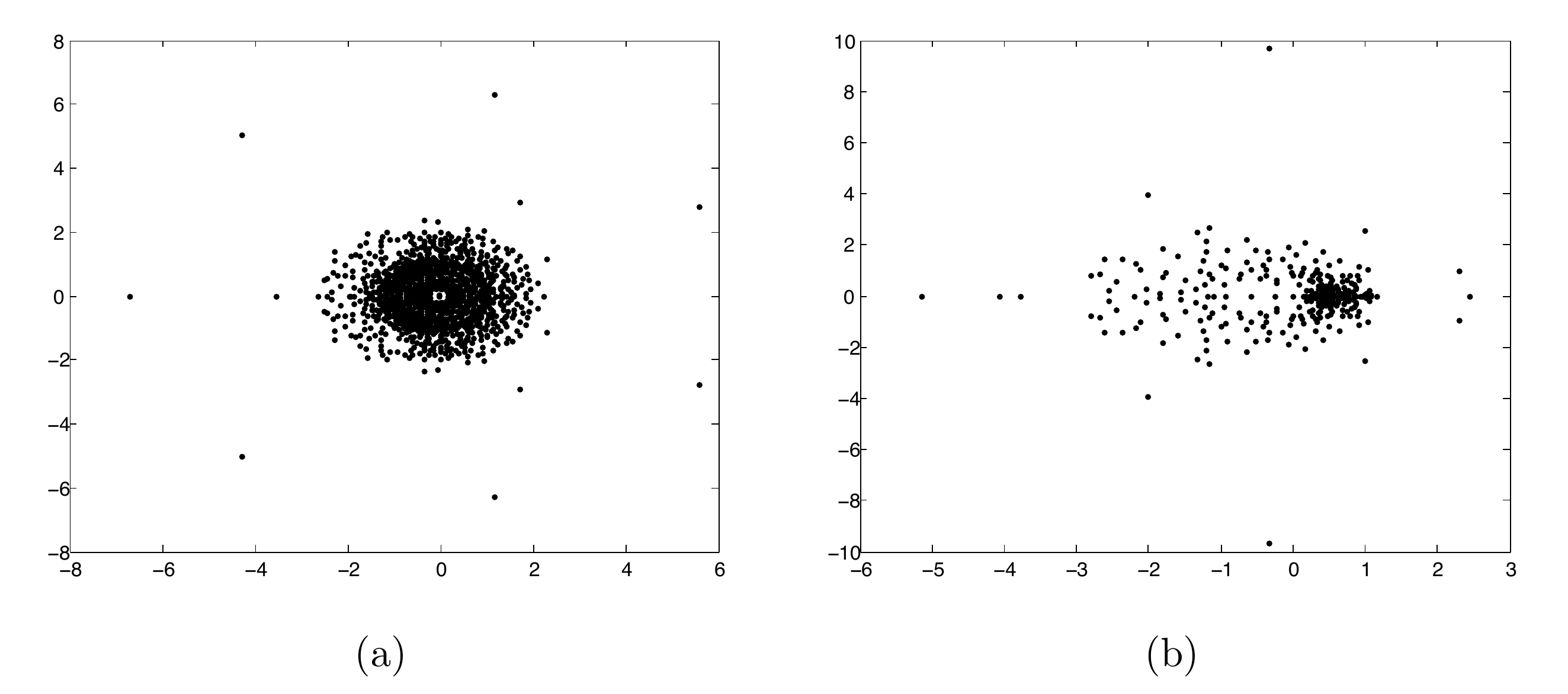}
\end{center}


\caption{(a) Eigenvalue distribution of the test matrix {\it mahindas}. (b) Eigenvalue distribution of the ILU($0$)-preconditioned {\it mahindas}. }
\label{fig:1-29-1}
\end{figure}

In this subsection, we demonstrate the effect of the deflation-subspace matrix $Z$ defined by (\ref{equ:11-24-1}) applied to the solution of
the following two test data downloaded from The University of Florida Sparse Matrix
Collection\footnote{http://www.cise.ufl.edu/research/sparse/matrices/}:
\begin{enumerate}
\item[(a)] {\it bcsstm27} from a mass matrix buckling problem. {\it bcsstm27}
is 
a $1224 \times 1224$ real symmetric and indefinite matrix $A$ with $56,126$ nonzero entries. As the right-hand side in (\ref{equ:11-19-1}), we set $b = A {\bf 1}$ where ${\bf 1} = [1, 1, \ldots, 1]^T$. A spectral plot for {\it bcsstm27} is in Figure \ref{fig:ellipse}(b).
\item[(b)] {\it mahindas} from an economic problem. {\it mahindas}
is a $1258 \times 1258$ real unsymmetric matrix $A$ with $7,682$ nonzero entries. Again, we set $b = A {\bf 1}$ as the right-hand side in (\ref{equ:11-19-1}). A spectral plot for {\it mahindas} is in Figure
\ref{fig:1-29-1}(a).
\end{enumerate}

All the computations were done in Matlab Version 7.1 on a Windows 7 machine with a Pentium 4 processor. An ILU preconditioner generated by the Matlab function $[L, U, P] = luinc(A, '0')$ was used for {\it mahindas}, namely, instead of solving (\ref{equ:11-19-1}), we solved
$$
\tilde{A} x = \tilde{b}
$$
where $\tilde{A} = L^{-1} PA U^{-1}$ and $\tilde{b} = L^{-1} P b$, and accordingly the $A$ and $b$ in (\ref{equ:7-31}) were replaced with $\tilde{A}$ and $\tilde{b}$ respectively. Since the $U$ factor obtained from {\it luinc} had some zeros along its main diagonal, we replaced those zeros by $1$ so that $U$ was invertible.
A spectral plot for $\tilde{A}$ is given in Figure
\ref{fig:1-29-1}(b).
On the other hand, we did not use any preconditioner for {\it bcsstm27}.

Numerical solutions with deflated restarted GMRES of the linear systems resulted from the discretization of the
two dimensional steady-state convection-diffusion equation
\begin{equation}\label{equ:11-24-10}
\begin{array}{rcl}
-[u_{xx} + u_{yy} + Re\, (p(x, y) u_x + q(x, y) u_y)] = f(x, y), & & (x, y) \in [0, 1]^2
\end{array}
\end{equation}
with Dirichlet boundary conditions were
studied in depth in \cite{dinkla}.
%
In \cite{dinkla}, two types of delation-subspace matrix $Z$ are used: eigenvectors obtained from the Matlab function $eig$, and algebraic subdomain vectors. The $Z$ of algebraic subdomain vectors works well for the fluid flow problem (\ref{equ:11-24-10}), but not for other problems. Accurately calculating eigenvalues of large linear systems, on the other hand, is very time-consuming. Therefore deflation with the $Z$ of true eigenvectors is not practicable.
Numerical experiments in \cite{dinkla} show that eigenvalues close to the origin hamper the convergence of a Krylov subspace method. Hence, deflation of these eigenvalues is very beneficial. Based on this observation, we chose in our experiments the $\Gamma$ in (\ref{equ:11-24-1}) to be a circle $D(c, r)$
with the center $c$ near the origin. For the $Y$ in (\ref{equ:11-24-1}),
we picked a random $Y \in {\mathbb R}^{N \times m}$ with $m$ not less than the exact number $s$ of eigenvalues inside $\Gamma$. We remark that an efficient stochastic estimation method of $s$ has been developed in \cite{futa}.
Moreover, we computed the integral in (\ref{equ:11-24-1}) by the Legendre-Gauss
quadrature
\begin{equation}\label{equ:11-24-13}
Z = \frac{r}{2} \int_{-1}^1 e^{\pi \theta \sqrt{-1}} ((c+r e^{\pi \theta \sqrt{-1}}) I - A)^{-1} Y d \theta
 \approx \frac{r}{2} \sum_{k = 1}^q \omega_k e^{\pi \theta_k \sqrt{-1}} ((c+r e^{\pi \theta_k \sqrt{-1}}) I - A)^{-1} Y,
\end{equation}
where $\omega_k$ and $\theta_k$ are the Legendre-Gauss weights and nodes on the interval $[-1, 1]$ with truncation order $q$.
In (\ref{equ:11-24-13}),
there are $mq$ linear systems $((c+r e^{\pi \theta_k \sqrt{-1}}) I - A) x = y_j$ to solve.
We solved each of them by BiCG with the stopping tolerance $tol = 10^{-10}$ and the maximum number of iterations $maxit = N$.


In our experiments, we performed the following three computations:
\begin{enumerate}
\item[\#1] Solve (\ref{equ:11-19-1}) without any deflation.
\item[\#2] Compute $Z$ through (\ref{equ:11-24-13}). Perform $QR$ factorization on $Z$: $Z = QR$ where $Q \in {\mathbb C}^{N \times m}$ and $R \in {\mathbb C}^{m \times m}$. Then
set $Z = Q$ which is the $Z$
in (\ref{equ:8-21-1}). Then solve (\ref{equ:7-31}).
\item[\#3] Use the Matlab function $eig$ to compute the eigenvectors $v_1, v_2, \ldots, v_s$ of $A$ whose associated eigenvalues lying inside $\Gamma$. Pick an $M \in {\mathbb R}^{s \times m}$ randomly, and set $Z = [v_1, v_2, \ldots, v_s] M$. Perform $QR$ factorization on $Z$: $Z = QR$ where $Q \in {\mathbb C}^{N \times m}$ and $R \in {\mathbb C}^{m \times m}$. Then
set $Z = Q$ which is the $Z$
in (\ref{equ:8-21-1}). Then solve (\ref{equ:7-31}).
\end{enumerate}

Due to that full GMRES is too expensive for us to use in terms of time and storage,
rather than use GMRES, we employed BiCG as the Krylov solver in the solution of (\ref{equ:11-19-1}) and (\ref{equ:7-31}). The initial guesses for BiCG were $x = 0$, and the stopping criteria were $\|b - Ax\|_2/\|b\|_2 < 10^{-7}$ for (\ref{equ:11-19-1})
and $\|Pb - PAx\|_2/\|Pb\|_2 < 10^{-7}$ for (\ref{equ:7-31}) respectively.

Numerical results are summarized in Table \ref{tab:11-24-1}.
In this table, the column titled with ``\#eig in $\Gamma$'' is a column of numbers of
eigenvalues of $A$ inside $\Gamma$.
The columns titled with ``\#iter'' are columns
of numbers of iterations, and the columns with ``Err''
are columns of true relative errors $\|b - Ax\|_2/\|b\|_2$ or $\|Pb - PAx\|_2/\|Pb\|_2$.

In these two experiments, BiCG essentially did not converge. With an appropriate eigenvalue-deflation, however, the situation was changed significantly.
The most expensive part in the proposed method is clearly the computation of the $Z$ in (\ref{equ:11-24-13}).
In next subsection, we describe the state-of-the-art parallel multigrid methods which may be applied to the computation.

\begin{table}
\centering
\footnotesize{
{\newcommand{\q}[1]{\mc{1}{|l|}{\small\tt #1}}
\noindent
\begin{tabular}{*{10}{|c}|} \hline
& Circle $\Gamma$ & & &
\mc{2}{|c|}{Computation \#1}& \mc{2}{|c|}{Computation \#2} & \mc{2}{|c|}{Computation \#3}
\\ \hline
Matrix & $(c, r)$ & \#eig in $\Gamma$ &$m$
& \#iter & Err  & \#iter & Err & \#iter & Err
\\ \hline\hline
 {\it bcsstm27}& $(0, 5)$ &$363$ &$400$
   &$1224000$  &$
   4.0 \times 10^{-6}$ & $763$ & $
   8.9 \times 10^{-8}$ & $266$& $
   9.8 \times 10^{-8}$
\\ \hline
{\it mahindas}&$(-1,1)$ &$31$ &
   $50$ &$1258000$  &$1.3$ & $3937$ & $5.0 \times 10^{-8}$ &$1555$ &$4.7 \times 10^{-8}$
\\ \hline
\end{tabular}}
}


\caption{A comparison of solving (\ref{equ:11-19-1}) and (\ref{equ:7-31}) by BiCG. For {\it mahindas}, a ILU($0$) preconditioner was applied. $\Gamma$ is a circle with center $c$ and radius $r$.
The $q$ in
(\ref{equ:11-24-13}) is $q = 2^7$.
}
\label{tab:11-24-1}
\end{table}


\subsection{Multigrid}
\label{Subsec:Problem:Multigrid}

We can formulate either geometric multigrid
\cite{GPAstrakhantsev_1971a,
NSBakhvalov_1966a,
ABrandt_1977b,
CCDouglas_1984a,
RPFedorenko_1961a,
RPFedorenko_1964a,
WHackbusch_1985a,
PWesseling_1992a}
or algebraic multigrid 
\cite{KStuben_2000a}
using the same
notation level to level using the abstract multigrid approach developed in
\cite{CCDouglas_1984c,
CCDouglas_1995a,
CCDouglas_JDouglas_DEFyfe_1994a,
REBank_CCDouglas_1985a,
CCDouglas_1994b,
CCDouglas_1995a,
CCDouglas_JDouglas_DEFyfe_1994a}.

Assuming the cost of the smoother (or rougher) on each level is $O(N_j)$,
$j=1,\cdots,k$, Algorithm MGC with $p$ recursions to solve problems on level
$k-1$ has complexity
\begin{equation}
\label{eqn:multigrid-complexity}
W_{MGC}(N_k) =
\left\{
\begin{array}{ll}
O(N_k)          & 1 \leq p \leq \sigma\\
O(N_k \log N_k) & p = \sigma\\
O(N_k^{\log p}) & p > \sigma.\\
\end{array}
\right.
\end{equation}
Under the right circumstances, multigrid is of optimal order as a solver.

Consider the example (\ref{equ:11-24-10}) in \S\ref{Subsec:Problem:NumerExI}.
A simple geometric multigrid approximation to (\ref{equ:11-24-10}) produces
a very good solution in $4$~V Cycles or $2$~W cycles using
the deflated GMRES as the rougher.
Each V or W Cycle is $O(N_k)$.
Hence, we have an optimal order solver for (\ref{equ:11-24-10}), which
would not be the case if we used BiCG or deflated GMRES on a single grid.

High performance computing versions of multigrid based on using
hardware acceleration with memory caches was extensively studied in the
early 2000's
\cite{CCDouglas_JHu_MKowarschik_URude_CWeiss_2000a}.

Parallelization of Algorithm MGC is straightforward
\cite{CCDouglas_1996c}.
\begin{itemize}
\item
For geometric multigrid, on each level~$j$,
data is split using a
domain decomposition paradigm.
Parallel smoothers (roughers) are used.
The convergence rate degrades from the standard serial theoretical rate,
but not by a lot, and scaling is good given sufficient data.
\item
For algebraic multigrid, the algorithms can be either straightforward
(e.g., Ruge-Studen \cite{JWRuge_KStuben_1985a}
or
Beck \cite{RBeck_1999a}) to quite complicated (e.g.,
AMGe \cite{GHaase_2000a}).
Solutions have existed for a number of years, so it is a matter of
choosing an exisiting implementation.
In some cases, using a tool like METIS or ParMETIS is sufficient to
create a domain decomposition-like system based on graph connections in
$A_j$, which reduces parallelization back to something similar to the
geometric case.
\end{itemize}
In many cases, the complexity of this type of parallel multigrid for $P$
processors becomes
\begin{equation}
\label{eqn:ParallelMGCcost}
W_{MGC,P}(N_k) = W_{MGC}(N_k)\log P / P.
\end{equation}


\section{Conclusions and Future Work}
\label{sect:future-work}

The novelties of this research include (i) we incorporate the
delation projector $P$ with the $Z$ described in (\ref{eq:projector}) and
(\ref{eq:gaussian}) into Krylov subspace methods to enhance the stability and
accelerate the convergence of the iterative methods for solving
ill-conditioned linear algebraic systems, and (ii) we will also implement robust and
efficient parallel multigrid methods for solving (\ref{eq:gaussian}) and
realize a software package for a wide variety of applications.

To our best knowledge, the constructions of most, if not all, deflation
subspace matrices $Z$ in the literature are problem dependent.
Further, some of them are ad-hoc, e.g., the algebraic subdomain deflation in
\cite{fv01}. The method proposed here is problem independent.


\section*{Acknowledgments}
\label{sect:acks}

This research was supported in part by National Science Foundation grants
ACI-1440610, ACI-1541392, and DMS-1413273.


%
\label{sect:bib}
 \bibliographystyle{plain}
\bibliography{iccs2016}






\end{document}